\documentclass[a4paper, 12pt]{article}
\usepackage{enumerate, theorem}
\usepackage{amsmath, amsfonts, amssymb}
\usepackage[height=22.5cm, width=15cm]{geometry}
\usepackage[all]{xy}
\usepackage{mathrsfs, yfonts}

\DeclareMathOperator{\C}{\mathbb{C}}

\newcommand{\parag}[1]{\paragraph{\sc{#1.}}}

\newtheorem{thm}{Theorem}[subsection]
\newtheorem{defn}[thm]{Definition}
\newtheorem{cor}[thm]{Corollary}
\newtheorem{prop}[thm]{Proposition}
\newtheorem{lemma}[thm]{Lemma}

\begin{document}

\title{Why (a,b)-modules, frescos and themes ?}

\author{Daniel Barlet\footnote{Barlet Daniel, Institut Elie Cartan UMR 7502  \newline
Universit\'e de Lorraine, CNRS, INRIA  et  Institut Universitaire de France, \newline
BP 239 - F - 54506 Vandoeuvre-l\`es-Nancy Cedex.France. \newline
e-mail : daniel.barlet@univ-lorraine.fr}.}

\maketitle

\parag{Abstract} This text is an advocacy for the use of (a,b)-modules (formal Brieskorn modules) in the study of the singularity of a holomorphic function on a complex manifold. It gives a short
an elementary overview on this simple algebraic tool.

\parag{AMS Classification} 32 S 05 - 34 E 05 - 32 S 30 - 32 S 40 - 34 A 30

\tableofcontents

\newpage

\section{Introduction}

I have studied the (a,b)-modules (formal Brieskorn modules)  for almost 30 years now and using this structure in many papers on singularity theory. Looking back to several of these papers I realized that may be I have not 
made enough efforts to convince colleagues that this structure is really interesting and important. The fact to inverse $\partial_z$ in the Weyl algebra $\C\langle z, \partial_z \rangle$ seems a priori
"anecdotal" and without serious reasons, mathematicians dont want to replace a "usual" differential equation by a structure using formal series in $\partial_z^{-1}$ and avoiding to use $\partial_z$.
The aim of the present text is to give some serious but simple arguments in favor of the use of this structure. I can summarize them in 3 points:
\begin{itemize}
\item In singularity theory this structure appears naturally and is also present on the space of asymptotic expansions of functions of the Nilsson class.
\item There is a very easy way to encode in this structure the regularity of the Gauss-Manin connection.  The notion of Bernstein polynomial of a regular (a,b)-module gives in this abstract setting  a definition of the logarithm of the eigenvalues of the monodromy.  This leads to the notion of "geometric (a,b)-module" which is defined by the regularity and the fact that the roots of the Bernstein polynomial are rational and negative. Moreover the structure of a regular (a,b)-module contains in a natural way an analog of a Hodge filtration (the $b-$adic filtration) and of a weight filtration associated to the nilpotent part of the monodromy (given by the semi-simple and co-semi-simple filtrations). Moreover we have nice functorial properties for geometric (a,b)-modules.
\item When we make a formal change of variable $\theta(z) \in \C[[z]]$ with $\theta(0) = 0$ and $\theta'(0)\not= 0$, which means that we replace "a"  by $\alpha := \theta(a)$ and  "b"  by $\beta := b.\theta'(a)$ we obtain a new (a,b)-module and this "change of variable" preserves many interesting properties like regularity, the Bernstein polynomial (and so geometric (a,b)-modules are stable by change of variable), semi-simplicity, monogenic 
$\tilde{\mathcal{A}}$ modules (frescos) ...  (see \cite{[B.12a]}).
\end{itemize}
This last point allows to use (a,b)-modules in the study of the singularities of a holomorphic map $f:  M \to C$ where $M$ is a complex manifold and $C$ is a smooth curve. To use an (a,b)-module structure at a critical value 
$p \in C$ of $f$, it is necessary to choose a local coordinate near $p$ on $C$ in order to have a complex valued function on $M$ near the singular fiber $f^{-1}(p)$.\\
 Then the existence and computation of (quasi-)invariants numbers by change of coordinate in the set of isomorphism classes of  geometric (a,b)-modules is an interesting challenge in such a context (see \cite{[B.12a]}).  \\
Let me conclude this introduction by saying that when the holomorphic function we are interested in depends on a holomorphic parameter, the consideration of a holomorphic family of geometric (a,b)-modules is helpful
for the study of period-integrals depending on such a parameter (see \cite{[B.15]}).\\

 \section{What is a (a,b)-module}
Consider a germ of holomorphic function $f : (\C^{n+1}, 0) \to (\C, 0)$ and also a germ of  $f-$relative holomorphic $(p+1)-$differential form $\omega \in \Omega^{p+1}_{\C^{n+1}, 0}$
satisfying $df\wedge \omega = 0$ and $d\omega = 0$.  For $p = n$ these two conditions are satisfied by any $\omega$ in  $ \Omega^{n+1}_{\C^{n+1}, 0}$.\\
 For each (compact) $p-$cycle  in the Milnor's fiber of $f$ we consider the (germ of)  {\bf period-integral}
$$ \varphi_\gamma(s) := \int_{\gamma_s} \omega/df $$
where $(\gamma_s)_s \in H$ is the horizontal family of $p-$cycles in the fiber of $f$ parametrized by the universal cover $H \to D^*$ of a small enough punctured disc $D^*$ with center $0$
in $\C$. This germ of multivalued function with finite determination has an asymptotic expansion at $s = 0$ of the type
$$ \sum_{r \in \Lambda, j \in [0, p], m \in \mathbb{N}} c^{r, j}_m(\gamma).s^{m+r}.(Log\, s)^j $$
where $\Lambda$ is a finite set in $\mathbb{Q}\cap ]-1, 0]$ which converges near the origin (in suitable sectors).
In fact, $(\gamma, s) \mapsto \varphi_\gamma(s)$ may be viewed as a germ of multivalued function on the product $H_p(F, \C)\times H $ where $F$ is the Milnor's fiber of $(f, 0)$. \\
For each given $p \in [1, n]$  this function is solution on $(\C, 0)$ of a "regular singular" differential system which is independent of the choices of $\omega$ and $\gamma$ which is the {\bf Gauss-Manin system} in degree $p$.
 The classical case of  a germ $(f, 0)$ with an isolated singularity at the origin\footnote{where only the degree $p = n$ is non trivial} shows that the "good point of view" to look at this differential system is the notion of {\bf Brieskorn module} which replaces the meromorphic connection (on the vector bundle $R^nf_*(\C)\otimes \mathcal{O}_{D^*}$)  by the operations  "product by $s$" and  "$\partial_s^{-1}$", that is to say the primitive in "$s$"  vanishing at $s = 0$.\\
 
 This leads to consider a new structure: the notion of "geometric (a,b)-module" corresponding to the formal completion of the Brieskorn module. The regularity of the Gauss-Manin connection implies that the asymptotic expansions we are interested in always converge, so there is no lost of generality to look at them as formal developments. The  "geometric" condition encodes in the structure of an  (a,b)-module both the regularity of the Gauss-Manin connection and the fact that the roots of  the Bernstein polynomial of $f$ are negative rational numbers.\\
 
 The notion of an (a,b)-module itself is very simple:\\
 Consider the  unitary $\C-$algebra $\mathcal{A} := \C < a, b >$ with two variables a and b  satisfying  the commutation relation $a.b - b.a = b^2$ corresponding to the relation $\partial_s.s - s.\partial_s = 1$ because "a" will act as the "multiplication by $s$" and "b" will act as the "primitive vanishing at $s=0$". Then complete $\mathcal{A}$ for the $b-$adic filtration of $\mathcal{A}$ defined by the two sided ideals $b^m.\mathcal{A} = \mathcal{A}.b^m, m \in \mathbb{N}$. We obtain the unitary $\C-$algebra $\tilde{\mathcal{A}}$ described as follows:
  $$ \tilde{\mathcal{A}} := \{ \sum_{m = 0}^{+\infty}  \quad b^m.P_m(a), \ {\rm with} \quad  P_m \in \C[a] \quad \forall m \in \mathbb{N} \}$$
 with the product defined by the relation $a.b - b.a = b^2$ which implies \begin{equation}
 a.S(b) = S(b).a + b^2.S'(b)
 \end{equation}
  for any $S \in \C[b]$, where $S'$ is the "usual" derivative of the polynomial $S$. Inductively, this relation allows to write
 any element in $ \tilde{\mathcal{A}}$ on the form:
  $$\sum_{m = 0}^{+\infty} \quad  P_m(a).b^m, \ {\rm with} \quad P_m \in \C[a] \quad \forall m \in \mathbb{N}$$
 
 Note that $\C[[b]]$ and $\C[a]$  are  commutative sub-algebras in $\tilde{\mathcal{A}}$.\\
 
 \begin{defn}\label{(a,b) 1}
 An {\bf (a,b)-module} is a $\tilde{\mathcal{A}}-$module which is {\bf free and finite rank} on the sub-algebra $\C[[b]]$.
 \end{defn}
 
 This kind of object is very easy to construct, thanks to the following elementary lemma.
 
 \begin{lemma}\label{(a,b) 2}
 Let $E := \oplus_{j = 1}^k \C[[b]].e_j$ be a free rank $k$ $\C[[b]]-$module. For any choice of \ $x_1, \dots, x_k$ in $E$ there exists an unique structure of (a,b)-module on $E$ such that
 \begin{enumerate}
\item The action of $\C[[b]]$ on $E$ corresponds to the initial structure of $\C[[b]]-$module of $E$.
 \item  For each $j \in [1, k]$ we have $a.e_j = x_j$ and the action of  "a" is continuous with respect to the $b-$adic filtration $(b^m.E)_{m \in \mathbb{N}}$ of $E$.
 \end{enumerate}
 \end{lemma}
 
 \parag{Proof} This is an easy consequence of the formula $(1)$ extended to any $S \in \C[[b]]$ thanks to the hypothesis that "a"  is continuous for the $b-$adic filtration of $E$.$\hfill \blacksquare$\\
 
 It is shown in \cite{[B.93]} proposition 1.1 that for any $A \in End_{\C}(\C^k)\otimes_{\C} \C[[z]]$ there exists a unique rank $k$ (a,b)-module $E := \oplus_{j=1}^k \C[[b]].e_j$ such that the basis $(e)$ satisfies
 $a.(e) = A(a).b.(e) $;  so this (a,b)-module  represents the simple pole formal differential system 
 $$ z.\frac{dF}{dz}(z) = A(z).F(z) .$$
 
\section{Geometric (a,b)-modules}

 First we shall define "regular".
 
 \begin{defn}\label{regular 1}
 We shall say that an (a,b)-module $E$ is  a {\bf simple pole} (a,b)-module when it satisfies $a.E \subset b.E$.\\
 We shall say that an  (a,b)-module $E$  is {\bf regular} when it is a sub-(a,b)-module of a simple pole (a,b)-module.
 \end{defn}
 
 The simplest way to find if an (a,b)-module $E$ is regular is to consider its saturation by $b^{-1}.a$. It is easy to see that 
 \begin{equation}
  \tilde{E} := \sum_{j=0}^\infty  (b^{-1}.a)^j.E \subset E \otimes_{\C[[b]]}\C[[b]][b^{-1}]
  \end{equation}
   is stable\footnote{the action of "a" is defined on $\C[[b]][b^{-1}]$  by the commutation relation $a.b^{-1} =  b^{-1}.a -1$.}  by $a$, so $\tilde{E}$  is a $\tilde{\mathcal{A}}-$module which has no $b-$torsion and satisfies $a.\tilde{E} \subset b.\tilde{E}$ because $a.(b^{-1}.a)^j = b.(b^{-1}.a)^{j+1}$. So when $\tilde{E}$ is a finite rank $\C[[b]]-$module, it is a simple pole (a,b)-module and $E$ is regular. Conversely, it is not difficult to see that $\tilde{E}$ is the smallest $\C[[b]]-$module containing $E$ and stable by $b^{-1}.a$. So $E$ is regular {\bf if and only} if $\tilde{E}$ has finite rank over $\C[[b]]$.\\
   
   \begin{defn}\label{Bernstein}
   Let $E$ be a simple pole (a,b)-module. Then the {\bf Bernstein polynomial} of $E$ is the {\bf minimal polynomial} of the $\C-$endomorphism $(-b^{-1}.a)$ acting on the finite dimensional vector space $E\big/b.E$.\\
   When $E$ is a regular (a,b)-module the {\bf Bernstein polynomial} of $E$ is, by definition, the Bernstein polynomial of $\tilde{E}$, where $\tilde{E}$ is the saturation  of $E$ by $b^{-1}.a$.
   \end{defn}
   
   Of course, this definition comes from the fact that when the germ $(f, 0)$ has an isolated singularity, the Bernstein polynomial of $f$ at $0$  is equal to $(\lambda + 1).B_E(\lambda)$ where $B_E$ is the Bernstein polynomial of the (a,b)-module $E$  which is the $b-$completion of the Brieskorn module of $(f,0)$. The factor $(\lambda+1)$ corresponds to the Bernstein polynomial of a non singular germ.
   
\begin{defn}\label{geometric}
A regular (a,b)-module  $E$ is called {\bf geometric} when its Bernstein polynomial has {\bf negative rational roots}.
\end{defn}

The rationality of the roots are the counterpart of the fact that the local monodromy of a germ $(f, 0)$ is quasi-unipotent. The negativity is the counterpart of Malgrange positivity result (see \cite{[Ma.74]}). Of course, both encode the famous theorem of Kashiwara \cite{[K.76]} which says that the roots of the Bernstein polynomial of an arbitrary germ $(f, 0)$ has negative rational roots.\\

The regularity of the Gauss-Manin connection has the following counter-part which is often useful.

\begin{prop}\label{a-complete}
Let $E$ be a regular (a,b)-module. Then $E$ is complete for the $a-$filtration  $(a^m.E)_{m \in \mathbb{N}}$.\\
So any geometric (a,b)-module is a module over the algebra of formal power series
$$\hat{\mathcal{A}} := \{ P = \sum_{p,q \geq 0} c_{p,q}.a^p.b^q \}$$
which contains $\C[[a]]$ and $\C[[b]]$ as commutative sub-algebras.
\end{prop}

Note that any $P \in \hat{\mathcal{A}} $ which has a non zero constant term (i.e. $c_{0,0} \not= 0$) is invertible in $\hat{\mathcal{A}} $. So its action on any geometric (a,b)-module $E$ is bijective.\\

The {\bf general  existence theorem}  of geometric (a,b)-modules\footnote{see \cite{[B.06]}  paragraph 2.3. The global case, when $f$ is proper, which associates  geometric (a,b)-modules to the hyper-cohomology of the complex $(\hat{K}er(df)^\bullet, d^\bullet)$ is given in  \cite{[B.13a]} theorem 4.3.4. See also  \cite{[B.12b]} theorem 2.3.1 for the relative case.}  shows that for each degree $p \in [0, n]$  the local Gauss-Manin connection of an arbitrary germ $(f, 0)$ is associated to a geometric (a,b)-module $E^p$.\\
This (a,b)-module $E^p$ is obtained as the quotient by its $b-$torsion of the germ of the $p-$th cohomology sheaf of the complex of sheaves $(\hat{K}er(df)^\bullet, d^\bullet)$ defined on the hyper-surface $\{ f = 0 \}$ in an ambiant complex manifold $M$  by defining $\hat{K}er(df)^m $ as the $f-$completion of the kernel of the sheaf map $(\wedge df) : \Omega_M^m \to \Omega_M^{m+1}$.\\

The next theorem explains that any geometric (a,b)-module can be realized as a sub-module of the asymptotic expansion (a,b)-module $\Xi^N_\Lambda$ defined as follows:\\

Let $\Lambda \subset \mathbb{Q}\cap ]0, 1]$ be a finite subset and let $N$ be a non negative integer. Then consider the free $\C[[b]]-$module of finite rank
$$ \Xi^N_\Lambda := \sum_{\lambda \in \Lambda, j \in [0, N]} \  \C[[b]].s^{\lambda-1}.\frac{(Log\, s)^j}{j!} $$
where the action of "a" is given by left multiplication by $s$ with the commutation relation given by the formula $(1)$. Interpreting "b" as the primitive without constant, it is easy to see that
$$ \Xi^N_\Lambda \simeq  \sum_{\lambda \in \Lambda, j \in [0, N]} \  \C[[s]].s^{\lambda-1}.\frac{(Log\, s)^j}{j!} $$
is the "usual" $\C[[s]]-$module of asymptotic expansions.\\
For a finite dimensional vector space $V$, we have a natural structure of $\tilde{\mathcal{A}}-$module on \  $ \Xi^N_\Lambda\otimes_{\C} V$  \ by letting $a$ and $b$  act as $a\otimes Id_V$ and $b \otimes Id_V$. Moreover, it is easy to see that this simple pole (a,b)-module is geometric for any choices of $\Lambda, N$ and $V$. 

\begin{thm}\label{D-A 1}
Any geometric (a,b)-module  $E$ may be realized as a sub-(a,b)-module of \ $\Xi^{N}_\Lambda\otimes_{\C} V$ \ for $N$ large enough and $\Lambda$ a finite subset in $\mathbb{Q}\cap ]0, 1]$, where $V$ is a complex vector space of dimension
at most equal to $rk(E)$.
\end{thm}

For a proof see \cite{[B.09b]}  theorem 4.2.1.

Let now $\Lambda $ be any subset in $\mathbb{Q}/\mathbb{Z}$. For any regular (a,b)-module $E$ let $E_\Lambda$ be the biggest sub-(a,b)-module $G$ in $E$ such that the roots of the Bernstein polynomial  $B_G$ of  $G$ are modulo $\mathbb{Z}$  in $-\Lambda$. then we have:

\begin{lemma}\label{primitif 1}
For any regular (a,b)-module $E$  and any  $\Lambda \subset \mathbb{Q}/\mathbb{Z}$, $E_\Lambda$ is a normal sub-(a,b)-module of $E$ and the Bernstein polynomial of the quotient $E/E_\Lambda$ has no root which is in $-\Lambda$ modulo $\mathbb{Z}$.
\end{lemma}

We shall call $E_\Lambda$ the {\bf $\Lambda-$primitive part} of $E$.\\
If $E$ is realized as a sub-(a,b)-module of $\Xi^{N}_M\otimes V$ then the $\Lambda-$primitive par of $E$ is given by the intersection
$$ E_{[\Lambda]} = E \cap \big( \Xi^{N}_{[\Lambda]}\otimes V\big).$$


\section{Semi-simplicity}

In order to study the structure of a regular (a,b)-module, we start with the classification of the rank 1 (a,b)-module. They are classified by a complex number. The rank 1 (a,b)-module corresponding to $\lambda \in \C$ is denoted $E_\lambda$ and is  defined by $E_\lambda := \C[[b]].e_\lambda$ where $a.e_\lambda = \lambda.b.e_\lambda$.\\
It is obvious that such an (a,b)-module has no non trivial normal submodule. In fact any non zero sub-(a,b)-module of $E_\lambda$ is equal to $b^m.E_\lambda$ which is isomorphic to $E_{\lambda+m}$ for some $m \in \mathbb{N}$. So we consider $E_\lambda, \lambda \in \C$ as the simple (a,b)-modules\footnote{There is an interesting analogy between regular (a,b)-modules and pairs $(V, f)$ of a finite dimensional complex vector space $V$ with an $\C-$endomorphism.}. Then we defined  a {\bf semi-simple} (necessarily regular) (a,b)-module $E$  as any  (a,b)-module isomorphic to a sub-module of a finite direct sum $\oplus_{j=1}^k  E_{\lambda_j}$ for some $\lambda_j \in \C$.\\
The following result is proved in \cite{[B.13b]} section 2.

\begin{prop}\label{ss-filt.}
In any regular (a,b)-module there exists a finite strictly increasing filtration $(S_i(E))_{i \in [0, d]}$ by normal sub-(a,b)-modules in $E = S_d(E)$ such that for each $i \in [1, d]$, $S_i(E)/S_{i-1}(E)$ is the biggest semi-simple sub-(a,b)-module inside $ E/S_{i-1}(E)$ (with the convention $S_0(E) := \{0\}$). It is called the {\bf semi-simple filtration} of $E$.\\
Dualy\footnote{We shall not describe here the notion of dual (a,b)-module for sake of brevity.} in any regular (a,b)-module there exists a finite strictly decreasing filtration $(\Sigma_h(E))_{h \in [1, \delta]}$ by normal sub-(a,b)-modules in $E := \Sigma_0(E)$ such that for each $h \in [1, \delta]$, $\Sigma_h(E)$ is the biggest normal  sub-(a,b)-module inside $\Sigma_{h-1}(E)$ such that the quotient $\Sigma_{h-1}(E)/\Sigma_h(E)$ is semi-simple. It is called the {\bf co-semi-simple filtration} of $E$.
\end{prop}

When a geometric (a,b)-module $E$ is embedded in some \ $\Xi^N_\Lambda \otimes V$ then the semi-simple filtration $S_j(E)$ is given by the intersection $E \cap (\Xi^j_\Lambda \otimes V)$. When the geometric (a,b)-module $E$  is a $\tilde{\mathcal{A}}-$sub-module generated by one element inside some $\Xi^j_\Lambda \otimes_{\C} V$  (so is a fresco, see the next section)  it is also easy to compute the co-semi-simple filtration $\Sigma_j(E)$ as the intersection of $E$ with $E \cap (\Xi^{N-j}_\Lambda \otimes V)$ assuming that $N$ is minimal for such an embedding of $E$.\\
In the analogy between regular (a,b)-modules and pairs $(V, f)$ where $V$ is a finite dimensional complex vector space and $f$ an endomorphism of $V$, the semi-simple filtration is the analog of the filtration $(Ker N^m)$ of $V$, where $N$ is the nilpotent part of $E$ (the analogy for the co-semi-simple filtration  uses duality).

\section{Why Frescos ?}

The Gauss-Manin system which corresponds to geometric (a,b)-modules controls the period-integrals associated to any (relative and closed) holomorphic differential form and any cycle $\gamma$ in the homology of the Milnor fiber of $(f, 0)$. But if we are interested in the period-integrals associated to a {\bf given class} of holomorphic $p-$form $\omega$ such that $df\wedge \omega = 0$ and $d\omega = 0$, we shall obtain much more precise informations if we may find a differential equation which has only these period-integrals as solutions. In term of (a,b)-module this corresponds to the fact that we want to consider only the sub-(a,b)-module $\tilde{\mathcal{A}}.[\omega] \subset E^p$ inside the geometric (a,b)-module $E^p$ associated to the Gauss-Manin connection of $(f, 0)$ in degree $p$ (see above). As a sub-(a,b)-module of a geometric (a,b)-module is again a geometric (a,b)-module, this leads to the following definition.

\begin{defn}\label{Fresco}
A geometric (a,b)-module which is generated as a $\tilde{\mathcal{A}}-$module by {\bf one} generator is called a {\bf fresco}.
\end{defn}

Of course, this is equivalent to the fact that $E\big/a.E + b.E$ is a  $1-$dimensional complex vector space, when $E$ is a geometric (a,b)-module.\\

This notion has several interesting properties:

\begin{lemma}\label{basis}
Let $F$ be a fresco with  generator $e$ over the algebra $\tilde{\mathcal{A}}$ and with $\C[[b]]-$rank $k$. Then $e, a.e, \dots, a^{k-1}.e$ is a $\C[[b]]-$basis of $F$.
\end{lemma}

\begin{lemma}\label{quotient and +}
Let $G \subset F$ be a normal\footnote{So for any $x \in F$,  $b.x \in G$ implies $x \in G$.} sub-(a,b)-module of a fresco $F$. Then $G$ is a fresco. Moreover,  $F/G$ is also a fresco.
\end{lemma}

For a proof see \cite{[B.13b]} lemma 1.1.4.

\begin{lemma}\label{Bernstein 2}
When $F$ is a fresco, the Bernstein polynomial of $F$ is equal to the {\bf characteristic polynomial} of $-b^{-1}.a$ acting on $\tilde{F}/b.\tilde{F}$ (where $\tilde{F}$ is the saturation of $F$ by $b^{-1}.a$).
\end{lemma}

The following  structure theorem is proved in \cite{[B.09b]} theorem 3.4.1. 

\begin{thm}\label{Fresco structure}
Let $F$ be a fresco with $\C[[b]]-$rank $k$. Then $F$ is isomorphic to a quotient $\tilde{\mathcal{A}}\big/\tilde{\mathcal{A}}.\Pi$ where 
$$ \Pi := (a-\lambda_1.b)S_1^{-1}.(a - \lambda_2.b).S_2^{-1}...S_{k-1}^{-1}.(a - \lambda_k.b) $$
where for each $j \in [1, k] , S_j$ is  in $\C[b]$ and   satisfies $S_j(0) = 1$, and where $\lambda_1, \dots, \lambda_k$ are rational numbers.
\end{thm}

Note that in such a presentation of the fresco $F$ the element
 $$P_F := (a - \lambda_1.b).(a - \lambda_2.b)...(a - \lambda_k.b)$$
 which is the initial form in (a,b) of $\Pi$ {\bf depends only on the isomorphism class} of  $F$ and is related to the Bernstein polynomial $B_F$ of $F$ by the relation
$$ P_F = b^{k}.B_F(-b^{-1}.a) $$
in the algebra $\mathcal{A}[b^{-1}]$ (see the theorem 3.2.1 in \cite{[B.09b]}).
We call the homogeneous element $P_F$  of degre $k := rk(F)$ in (a,b) the {\bf Bernstein element} of $F$. Then note that the roots of $B_F$ are the (rational negative) numbers $-\lambda_j + k - j$. 
So  $B_F(x) = \prod_{j=1}^k (x + \lambda_j +j - k)$. This implies that we have $\lambda_j > k-j ,\ \forall j \in [1, k]$.\\

The following corollary of the previous results shows that, in the case of a fresco, the Bernstein polynomial is  easier to compute than for a general geometric (a,b)-module.

\begin{cor}\label{suite exacte}
Let $0 \to G \to F \to H \to 0$ be an exact sequence of frescos.  Then we have the equality  $P_F = P_G.P_H$ in the algebra $\mathcal{A}$ between Bernstein elements. This gives the relation
$$  B_F(x) = B_G(x - rk(H)).B_H(x) \quad {\rm in} \ \C[x]$$
between the Bernstein polynomials.
\end{cor}

Another important result on frescos is given by the following theorem (see \cite{[B.13b]} theorem 1.2.5)

\begin{thm}\label{Principale J.H}
Let $F$ be a fresco. Then $F$ admits an {\bf unique Jordan-H\"older sequence}
 $$0 = F_0 \subset F_1 \subset F_2 \subset \dots \subset F_k = F $$
 where, for each $j \in [0,k]$, $F_j$ is a normal fresco of rank $j$ such that $F_{j}/F_{j -1}\simeq E_{\lambda_{j}}$ {\bf such that the sequence $(\lambda_j +j)_{j \in [1, k]}$ is non decreasing}. \\
 It is called the {\bf principal Jordan-H\"older sequence} of the fresco\footnote{Recall that $E_\lambda$ the rank 1 regular  (a,b)-module defined by $E_\lambda := \C[[b]].e_\lambda$ where $a.e_\lambda = \lambda.b.e_\lambda$.}.
 \end{thm}
 
\parag{Remarks}  
\begin{enumerate}
\item Any regular (a,b)-module admits a Jordan-H\"older sequence, but, even with the condition that the sequence $(\lambda_j + j)_{j \in [1, k]}$ is non decreasing, it is not unique. in general.
\item For a fresco $F$ to admit a Jordan-H\"older sequence with $\lambda_j+j$ strictly decreasing is equivalent to semi-simplicity.
\item A $\lambda-$primitive theme (see below) has a unique Jordan-H\"older sequence. So it is its principal Jordan-H\"older sequence. This shows that a non zero semi-simple $\lambda-$primitive theme has rank $1$.
\end{enumerate}

In term of asymptotic expansions, the notion of fresco has the following characterization which is an easy corollary of the theorem \ref{D-A 1}.

\begin{cor}\label{D-A 2}
Any fresco $F$ is isomorphic to a $\tilde{\mathcal{A}}-$monogenic sub-module of some $\Xi^N_\Lambda\otimes V$ for some finite subset $\Lambda \subset \mathbb{Q} \cap[0, 1[$, some $ N \in \mathbb{N}$ and some finite dimensional complex vector space $V$.
\end{cor}

For the fresco $F := \tilde{\mathcal{A}}.[\omega] \subset E^p$ associated to the germ $(f, 0)$ and $\omega \in \Omega_{\C^{n+1}, 0}^p$ such that $df\wedge \omega = 0$ and $d\omega = 0$, we may choose $V := H^p(F, \C)$ where $F$ is the Milnor fiber of $(f, 0)$ and for each $\gamma \in H_p(F, \C) \simeq H^p(F, \C)^*$  the image by $Id\otimes \gamma$ of $[\omega]$  in $\Xi^N_\Lambda$ is given by the asymptotic expansion at $s= 0$  of the integral period $\varphi_\gamma$. This leads to the notion of theme.

\section{Themes}

When we are interested by a period-integral $\varphi_\gamma(s) := \int_{\gamma_s} \omega/df$ with a {\bf given class} of $d-$closed $f-$relative holomorphic $(p+1)-$differential form $\omega$ and a {\bf given class}  $\gamma \in H_p(F, \C)$, the asymptotic expansion at $s = 0$ of $\varphi_\gamma(s)$ defines a monogenic sub$-\tilde{\mathcal{A}}-$module of some $\Xi^N_\Lambda$. This will be called the {\bf theme} of this period-integral.\\

Among the frescos, themes are characterized by the fact that they may be realized in the asymptotic expansions with a $1-$dimensional vector space $V$. This is equivalent to ask that, for each eigenvalue of the monodromy\footnote{the monodromy may be defined directly on $\Xi^N_\Lambda$ by $T(\psi(s)) = \psi(s.e^{2i\pi})$ on the formal power series in $s = \exp(2i\pi.\sigma), \sigma \in H$; so $T$ is the translation by $2i\pi$ on $\sigma := Log\, s$.}, there is exactly one Jordan bloc. Fixing $\gamma \in H_p(F, \C)$ in the case of a period-integral, this Jordan bloc is generated by the component of $\gamma$ on the corresponding spectral subspace of  the monodromy  of $f$ acting on the $p-$th homology group $H_p(F, \C)$ of the Milnor fiber of $f$.\\

The following result is an "abstract version" of our previous remark on the monodromy of a period-integral associated to a given $\omega$ and a given $\gamma$ which also gives a characterization of $\lambda-$primitive themes.

\begin{thm}\label{primitif 2}
Let  $\lambda$ be in $\mathbb{Q}/\mathbb{Z}$ and let $\Theta$ be a $\lambda-$primitive theme. Then $\Theta$ has a unique normal rank $1$ sub-module. Conversely, any geometric (a,b)-module which has an unique rank $1$ normal sub-module is a $\mu-$primitive theme for some $\mu \in \mathbb{Q}/\mathbb{Z}$.
\end{thm}

For a proof see \cite{[B.14]} theorem 2.1.6.\\

As the quotient of a $\lambda-$primitive theme by a normal sub-(a,b)-module is again a $\lambda-$primitive theme the following corollary is immediate.

\begin{cor}\label{unique J-H}
A  $\lambda-$primitive theme has a unique Jordan-H\"older sequence. Conversely, any geometric (a,b)-module having a unique Jordan-H\"older sequence is a $\mu-$primitive theme for some $\mu \in \mathbb{Q}/\mathbb{Z}$.
\end{cor}

Of course, this J-H. sequence must be the principal one, so the sequence $\lambda_j + j$ is non decreasing.

\parag{Notation} Let $\Theta$ be a $\lambda-$primitive theme. Then in its Bernstein element 
$$P_\Theta := (a -\lambda_1)\dots (a -\lambda_k)$$
 we have $\lambda_{j+1} = \lambda_j + p_j - 1$ for each $j \in [1, k-1]$ with $p_j \in \mathbb{N}$ and $\lambda_1$ which is is in $\lambda + \mathbb{Z}$, must satisfy $\lambda_1  > k-1$. 
 So $P_\Theta$ is defined by $\lambda_1$ and the integers $p_1, \dots, p_{k-1}$. \\
We shall say that $\lambda_1, p_1, \dots, p_{k-1}$ are the {\bf fundamental data} of the  $\lambda-$primitive theme $\Theta$.\\
Defined for each $j \in [1, k-1]$ the vector space $V_j \subset \C[[b]]$ as follows:
\begin{itemize}
\item When $p_j + \dots + p_{k-1} \geq k-j$ let $q_j := p_j+ \dots + p_{j+h} $ where $h$ is the smallest integer such that $p_j+ \dots + p_{j+h} \geq k-j$. Then put
 $$V_j := \oplus_{i= 0}^{k-j-1} \  \C.b^{i} \oplus \C.b^{q_j} .$$
 \item When $p_j + \dots + p_{k-1} <  k-j$ put 
 $$ V_j := \oplus_{i= 0}^{k-j-1} \  \C.b^{i} .$$
 \end{itemize}
 Note that  $b^0 = 1$ and $b^{p_j}$ are always in $V_j$.
 Define now for each $j \in [1, k-1]$ the open set $W_j$ in $V_j$ given by the conditions $S \in V_j$ belongs to $W_j$ if and only if $S(0) = 1$ and the coefficient of $S$ on $b^{p_j}$ is not zero.\\
The theorem \ref{Fresco structure} may be specified as follows (see \cite{[B.14]} theorem 3.2.6).

\begin{thm}\label{canonical form}{\bf [Canonical form for $\lambda-$primitive themes]}\\
Let $\Theta$ be a $\lambda-$primitive theme with fundamental data $\lambda_1, p_1, \dots, p_{k-1}$. Then there exists $S_j \in W_j$ such that $\Theta$ is isomorphic to $\tilde{\mathcal{A}}\big/\tilde{\mathcal{A}}.\Pi$
where 
$$\Pi := (a- \lambda_1.b).S_1^{-1}.(a - \lambda_2.b) \dots  S_{k-1}^{-1}.(a - \lambda_k.b) .$$
Conversely, for any such $\Pi \in \tilde{\mathcal{A}}$, the quotient $\tilde{\mathcal{A}}\big/\tilde{\mathcal{A}}.\Pi$ is a $\lambda-$primitive theme.
\end{thm}

Of course, this allows to produce for each given fundamental data a  versal family of $\lambda-$primitive theme (see \cite{[B.15]}).\\

By definition a theme is a sub-(a,b)-module of some $\Xi^N_\Lambda$. So it is natural to ask when such a realization is unique. This leads to the notion of {\bf invariant theme}\footnote{in French : th\`{e}me stable.}. 

The following result gives two characterizations of invariant themes (see proposition 4.2.4 of \cite{[B.14]}):

\begin{thm}\label{invariant} 
Let $\Theta$ be a theme with rank $k$. The following properties are equivalent:
\begin{enumerate}
\item There exists a realization of $\Theta$ in some $\Xi^N_\Lambda$ which is invariant by the monodromy (see a previous footnote for a definition of the monodromy on the asymptotic expansion $\Xi^N_\Lambda$$)$.
\item There exists an {\bf unique} sub-theme of \ $\Xi := \cup_{\Lambda, N} \  \Xi^N_\Lambda$ which is isomorphic to $\Theta$.
\item The $\C-$vector space $Hom_{\tilde{\mathcal{A}}}(\Theta, \Theta)$ has dimension $k$.
\end{enumerate}
\end{thm}

The following  theorem for $[\lambda]-$primitive invariant theme is proved in \cite{[B.14]}.

\begin{thm}\label{unique}
Any invariant $[\lambda]-$primitive theme has an unique canonical form.
\end{thm}

This result implies that in the versal family of themes with given fundamental data (parametrized by the affine open set  $\prod_{j=1}^k \ W_j$; see the theorem \ref{canonical form}) an invariant theme appears only once.\\
 This leads to construct some {\bf universal} families of $\lambda-$primitive themes (see \cite{[B.15]}).\\
 
 For the notion of {\bf change of variable} in an (a,b)-module the reader may consult \cite{[B.12a]}.\\
 
 For { \bf concrete non trivial examples} of the use of frescos to analyse asymptotic expansions of period-integrals see \cite{[B.21]}.

\end{document}